\documentclass[11pt]{amsart}
\usepackage[margin=1.1in]{geometry}
\usepackage{amsmath,amssymb,amsthm}
\usepackage[hidelinks]{hyperref}

\theoremstyle{plain}
\newtheorem{lemma}{Lemma}
\theoremstyle{definition}
\newtheorem{claim}{Claim}
\newtheorem*{conjecture}{Conjecture}

\newcommand{\Irr}{\mathrm{Irr}}
\newcommand{\F}{\mathbb{F}}
\newcommand{\Tr}{\operatorname{Tr}}

\title{A counterexample to Conjecture C}
\author{Youmbai Riad}
\address{Department of Mathematics\\ Faculty of Exact Sciences\\ University of El Oued\\ PO Box 789, 39000 Echott\\ El Oued\\ Algeria}
\email{youmbai-riad@univ-eloued.dz}
\date{}
\begin{document}
\maketitle

\begin{abstract}
We describe a solvable group $G$ of order about $2.4\times10^{34}$ and a primitive irreducible character $\chi$ of $G$ such that $\chi(1)$ divides no conjugacy class size of $G$. This contradicts Conjecture~C of Isaacs, Keller, Meierfrankenfeld and Moret\'o and Conjecture~8.1 of Navarro's survey.
\end{abstract}

\paragraph{AI disclosure.}
The construction and proofs in this note were generated by Claude, an AI language model developed by Anthropic. They have not been independently verified.

\section*{ \centering Introduction}

\begin{conjecture}[{\cite{IKMM}, Conjecture~C; \cite{Navarro}, Conjecture~8.1}]
If $\chi$ is a primitive irreducible character of a finite group $G$, then $\chi(1)$ divides $|g^G|$ for some $g\in G$.
\end{conjecture}

Here $g^G$ is the conjugacy class of $g$ in $G$. The claim of this note is that the conjecture fails already for solvable groups.

\section{The construction}

\subsection{Primes}

Let $(p_1,\dots,p_6)=(3,7,13,17,23,37)$. Each $p_j$ is odd and $p_j\equiv 2$ or $3\pmod 5$. Hence $p_j$ has order $4$ modulo $5$, and $5$ divides $p_j^2+1$.

\subsection{An automorphism of order 5}

\begin{lemma}\label{lem:aut}
Let $p$ be an odd prime with $p\equiv 2$ or $3\pmod 5$. Then there exist an extraspecial group $E$ of order $p^5$ and exponent $p$, and an automorphism $\alpha$ of $E$ of order $5$, such that $\alpha$ centralizes $Z(E)$ and acts irreducibly on $W=E/Z(E)$.
\end{lemma}

\begin{proof}
Let $F=\F_{p^4}$ and $F_0=\F_{p^2}$, and write $\bar x=x^{p^2}$ for $x\in F$. Then $x\mapsto\bar x$ is an automorphism of $F$ of order $2$ with fixed field $F_0$. Let $\Tr\colon F\to\F_p$ be the trace, so that $\Tr(x)=x+x^p+x^{p^2}+x^{p^3}$. Then $\Tr(\bar x)=\Tr(x)$ for all $x\in F$. Also $\Tr$ is not identically zero, since it is given by a polynomial of degree $p^3<|F|$, and hence $\Tr$ maps $F$ onto $\F_p$.

We construct a symplectic form on $F$. Choose $\varepsilon\in F\setminus F_0$ and put $\gamma=\varepsilon-\bar\varepsilon$. Then $\gamma\ne0$ and $\bar\gamma=-\gamma$. Define
\[
B(u,v)=\Tr(\gamma u\bar v)\qquad(u,v\in F).
\]
Then $B$ is $\F_p$-bilinear. For $u,v\in F$ we have
\[
B(v,u)=\Tr(\gamma v\bar u)=\Tr\bigl(\overline{\gamma v\bar u}\bigr)=\Tr(\bar\gamma\bar v u)=-\Tr(\gamma u\bar v)=-B(u,v).
\]
Hence $2B(u,u)=0$, and so $B(u,u)=0$ because $p$ is odd. If $u\ne0$, then $v\mapsto\gamma u\bar v$ is a bijection of $F$, and since $\Tr$ is onto, there exists $v$ with $B(u,v)\ne0$. Thus $B$ is nondegenerate, and $(F,B)$ is a $4$-dimensional symplectic space over $\F_p$.

Let $T=\{t\in F^\times: t\bar t=1\}$, which is the set of $t$ with $t^{p^2+1}=1$. Since $F^\times$ is cyclic of order $p^4-1=(p^2-1)(p^2+1)$, the group $T$ is cyclic of order $p^2+1$. For $t\in T$ we have
\[
B(tu,tv)=\Tr(\gamma\, tu\,\bar t\bar v)=\Tr(t\bar t\,\gamma u\bar v)=B(u,v),
\]
so multiplication by $t$ preserves $B$. Now $p\equiv\pm2\pmod 5$, so $p^2\equiv-1\pmod 5$ and $5$ divides $p^2+1$. Therefore $T$ contains an element $\sigma$ of order $5$.

We claim that $\sigma$ acts irreducibly on $F$ (as an $\F_p$-space) by multiplication. Let $U\ne0$ be a $\sigma$-invariant subspace. Then $U$ is invariant under the subring $\F_p[\sigma]$ of $F$. A finite subring of a field is a field, so $\F_p[\sigma]=\F_{p^k}$ for some $k$ dividing $4$. Since $\sigma$ has order $5$, we have $\sigma\in\F_{p^k}$ if and only if $5\mid p^k-1$. As $p\equiv2,3\pmod5$, the order of $p$ modulo $5$ is $4$, so $k=4$ and $\F_p[\sigma]=F$. For $0\ne u\in U$ we then get $U\supseteq Fu=F$, so $U=F$.

Now let $W=F$ and define $E=W\times\F_p$ with the multiplication
\[
(u,s)(v,t)=\bigl(u+v,\ s+t+\tfrac12B(u,v)\bigr),
\]
where $\tfrac12$ is the inverse of $2$ in $\F_p$. Using the bilinearity of $B$, one checks directly that both $\bigl((u,s)(v,t)\bigr)(w,r)$ and $(u,s)\bigl((v,t)(w,r)\bigr)$ equal
\[
\bigl(u+v+w,\ s+t+r+\tfrac12B(u,v)+\tfrac12B(u,w)+\tfrac12B(v,w)\bigr),
\]
so the multiplication is associative. The element $(0,0)$ is the identity, and $(u,s)^{-1}=(-u,-s)$ because $B(u,u)=0$. Thus $E$ is a group of order $p^5$. The set $Z=\{0\}\times\F_p$ is central, and $(u,s)\mapsto u$ is a homomorphism of $E$ onto $W$ with kernel $Z$, so $E/Z\cong W$. Since $B$ is skew, we get $(u,s)(v,t)=(v,t)(u,s)\,(0,B(u,v))$, so
\[
[(u,s),(v,t)]=(0,B(u,v)).
\]
As $B$ is nondegenerate, $(u,s)$ is central if and only if $u=0$, so $Z(E)=Z$. Also $B\ne0$, so $E'=Z$. Hence $E$ is extraspecial of order $p^5$. Finally, since $B(u,u)=0$, induction on $k$ gives $(u,s)^k=(ku,ks)$, so $(u,s)^p=(0,0)$ and $E$ has exponent $p$.

Define $\alpha\colon E\to E$ by $\alpha(u,s)=(\sigma u,s)$. This is a bijection, and since $B(\sigma u,\sigma v)=B(u,v)$,
\[
\alpha\bigl((u,s)(v,t)\bigr)=\bigl(\sigma u+\sigma v,\ s+t+\tfrac12B(u,v)\bigr)=\alpha(u,s)\,\alpha(v,t).
\]
So $\alpha$ is an automorphism. It fixes $Z$ pointwise, and $\alpha^5=1$ because $\sigma^5=1$, while $\alpha\ne1$ because $\sigma\ne1$. Hence $\alpha$ has order $5$. Under the identification $E/Z\cong W=F$, $\alpha$ acts as multiplication by $\sigma$, which is irreducible by the claim above.
\end{proof}

For each $j$, let $E_j$ be the group $E$ given by Lemma~\ref{lem:aut} for $p=p_j$, and let $\alpha_j$ be the corresponding automorphism of $E_j$. Put $Z_j=Z(E_j)$, of order $p_j$, and $W_j=E_j/Z_j\cong\F_{p_j}^4$, with the nondegenerate alternating form induced by commutators.

\subsection{The acting group and the group $G$}

Let $A=C_5\times C_5$. It has exactly six subgroups of order $5$; label them $K_1,\dots,K_6$. Let $A$ act on $E_j$ through the quotient $A/K_j\cong C_5$, a generator acting as the automorphism $\alpha_j$. Thus:
\begin{itemize}
\item $K_j$ acts trivially on $E_j$;
\item every $a\in A\setminus K_j$ acts on $W_j$ without nonzero fixed points;
\item $A$ acts trivially on $Z_j$.
\end{itemize}
To see the second item, let $a\in A\setminus K_j$. Then $a$ acts on $W_j$ as $\alpha_j^k$ with $k\not\equiv0\pmod5$, and $\langle\alpha_j^k\rangle=\langle\alpha_j\rangle$ because $5$ is prime. So the fixed space of $a$ on $W_j$ is $\langle\alpha_j\rangle$-invariant, hence it is $0$ or $W_j$ by Lemma~\ref{lem:aut}. It is not $W_j$, because $\alpha_j$ acts nontrivially on $W_j$: a trivial action would leave every line invariant, contradicting irreducibility since $\dim W_j=4$. Hence the fixed space is $0$.

Put
\[
N=E_1\times\cdots\times E_6,\qquad G=N\rtimes A,
\]
with $A$ acting componentwise. Then $|G|=25\prod_j p_j^5=24024018365307903996616487269711275\approx2.4\times10^{34}$, and $G$ is solvable. Since $5\nmid p_j$, the Sylow $p_j$-subgroup of $G$ is $E_j$, and it is normal in $G$.

\subsection{The character}

For each $j$ fix a faithful linear character $\lambda_j$ of $Z_j$, and let $\theta_j\in\Irr(E_j)$ be the unique irreducible character lying over $\lambda_j$; it has degree $p_j^2$. Since $A$ acts trivially on $Z_j$, it fixes $\lambda_j$, so for $a\in A$ the character $\theta_j^a$ also lies over $\lambda_j$, and hence $\theta_j^a=\theta_j$. Put
\[
\psi=\theta_1\times\cdots\times\theta_6\in\Irr(N),\qquad D=\psi(1)=\prod_j p_j^2=15598479159081.
\]
Then $\psi$ is $A$-invariant. Since $\gcd(|N|,|A|)=1$, $\psi$ extends to a character $\chi\in\Irr(G)$ (Isaacs \cite{Isaacs}, Cor.~8.16). So $\chi_N=\psi$ is irreducible and $\chi(1)=D$.

\section{Two claims}

\begin{claim}\label{cl:prim}
$\chi$ is primitive.
\end{claim}

\begin{proof}
Suppose $\chi=\eta^G$ with $\eta\in\Irr(H)$ and $H<G$.

(i) $|G:H|$ divides $\chi(1)=D$, and $|G:HN|$ divides $|G:N|=25$. Since $\gcd(D,25)=1$, we get $HN=G$.

(ii) Put $M=H\cap N$. As $G=HN$ there is a single double coset, so Mackey's theorem gives $\psi=(\eta_M)^N$. Since $\psi$ is irreducible, $\eta_M$ is irreducible and $\psi(1)=|N:M|\,\eta(1)$.

(iii) $Z(N)\le M$: for $z\in Z(N)$ we have $|\psi(z)|=\psi(1)\ne0$, while an induced character $\varphi^N$ vanishes outside $\bigcup_{x\in N}M^x$; a central element lying in $M^x$ lies in $M$.

(iv) Put $L=M/Z(N)\le W:=N/Z(N)=\bigoplus_jW_j$. Since $N'=\prod_jZ_j=Z(N)$, $N$ acts trivially on $W$ by conjugation, so $G$ acts on $W$ through $G/N\cong A$. As $G=HN$, $H$ maps onto $A$, and $L$ is $H$-invariant, hence $A$-invariant.

(v) The $W_j$ are $p_j$-groups for distinct primes, so $L=\bigoplus_j(L\cap W_j)$ with each $L\cap W_j$ an $A$-invariant subgroup of $W_j$. Since $A$ acts irreducibly on $W_j$, each $L\cap W_j$ is $0$ or $W_j$. So $L=\bigoplus_{j\in J}W_j$ for some set $J$.

(vi) Then $|N:M|=\prod_{j\notin J}p_j^4$ divides $\psi(1)=D$. The $p_j$-part of $D$ is only $p_j^2$, so $J$ contains every $j$. Hence $M=N$, so $N\le H$, and $H=HN=G$, a contradiction.

Hence $\chi$ is not induced from any proper subgroup.
\end{proof}

\begin{claim}\label{cl:class}
$\chi(1)=D$ divides no class size $|g^G|$.
\end{claim}

\begin{proof}
Write $g=na$ with $n=\prod_in_i\in N$ ($n_i\in E_i$) and $a\in A$.

(i) Every $a\in A$ lies in some $K_j$: $a=1$ lies in all of them, and each $a\ne1$ generates exactly one of the six subgroups of order $5$. Fix such a $j$.

(ii) $K_j$ acts trivially on $E_j$, so for $x\in E_j$ we get $x^g=(x^n)^a=x^{n_j}$. Hence $C_{E_j}(g)=C_{E_j}(n_j)$, which has index $1$ or $p_j$ in $E_j$ (in an extraspecial group every class has size $1$ or $p_j$).

(iii) $E_j$ is a normal Sylow $p_j$-subgroup of $G$, so every $p_j$-subgroup of $C_G(g)$ lies in $E_j$. Thus $C_{E_j}(g)=C_G(g)\cap E_j$ is a Sylow $p_j$-subgroup of $C_G(g)$, and
\[
|g^G|_{p_j}=|G:C_G(g)|_{p_j}=|E_j:C_{E_j}(g)|\le p_j .
\]

(iv) The $p_j$-part of $\chi(1)=D$ is $p_j^2>p_j$, so $D\nmid|g^G|$.
\end{proof}

Together, Claims~\ref{cl:prim} and~\ref{cl:class} say that $G$ is a solvable group with a primitive irreducible character whose degree divides no class size.

\end{document}